\documentclass[12pt] {article}
\usepackage{amsmath,amsthm}
\usepackage{amssymb,latexsym}
\title{The multiplicative structure on continuous polynomial valuations.}
\date{}
\author{ Semyon Alesker
\\  { \normalsize Department of Mathematics, Tel Aviv University, Ramat Aviv}
 \\  { \normalsize 69978 Tel Aviv,
Israel }
\\ {\normalsize e-mail: semyon@post.tau.ac.il}}

\def\RR{\mathbb{R}}
\def\CC{\mathbb{C}}

\def\ZZ{\mathbb{Z}}
\def\HH{\mathbb{H}}

\def\PP{\mathbb{P}}

\def\eps{\varepsilon}

\def\lam{\lambda}

\def\str{\longrightarrow}

\def\qed { Q.E.D. }
\def\pr{\partial}

\swapnumbers
\newtheorem{theorem}{Theorem}[section]
\newtheorem{corollary}[theorem]{Corollary}
\newtheorem{lemma}[theorem]{Lemma}
\newtheorem{proposition}[theorem]{Proposition}
\newtheorem{claim}[theorem]{Claim}
\theoremstyle{definition}
\newtheorem{example}[theorem]{Example}
\newtheorem{definition}[theorem]{Definition}
\newtheorem{remark}[theorem]{Remark}
\theoremstyle{proposition-definition}
\newtheorem{proposition-definition}[theorem]{Proposition-Definition}

\def\cf{{\cal F}}

 \def\ce{{\cal E}} \def\cf{{\cal F}}
\def\cg{{\cal G}}  
 \def\ck{{\cal K}}

\def\pvd{PVal_d(V)}
\def\ond{\Omega^{n}_{d}}
\def\pl{\PP_+(V^*)}

\def\gr{{}\!^ {\textbf{R}} Gr}

\begin{document}
\maketitle \setcounter{section}{-1}
\begin{abstract}
A canonical structure of a commutative associative filtered
algebra with unit is introduced on the space of polynomial smooth
valuations, and its properties are studied. Induced structure on
the subalgebra of translation invariant smooth valuations has
especially nice properties (it has the structure of the Frobenius
algebra). We also present some applications.
\end{abstract}
\section{Introduction.} The purpose of this paper is to introduce
a canonical structure of a commutative associative filtered
algebra with unit on the space of polynomial smooth valuations and
study its properties. The induced structure on the subalgebra of
translation invariant smooth valuations has especially nice
properties (it is the structure of the Frobenius algebra). We will
also present some applications. Part of the results of this paper
was announced in \cite{alesker-icm}.

Let us describe our main results in more detail. Let us recall
first of all some notation and definitions. Let $V$ be a real
vector space of finite dimension $n$. Let $\ck(V)$ denote the
class of all convex compact subsets of $V$. If we fix on $V$ a
Euclidean metric then we can define the Hausdorff metric $d_H$ on
$\ck(V)$ as follows: $$d_H(A,B):=\inf\{ \eps >0|A\subset (B)_\eps
\mbox{ and } B\subset (A)_\eps\},$$ where $(U)_\eps$ denotes the
$\eps$-neighborhood of a set $U$. However the topology on $\ck(V)$
does not depend on the choice of a Euclidean metric on $V$.
Moreover, equipped with this topology, $\ck(V)$ becomes a locally
compact topological space (Blaschke's selection theorem).

\begin{definition}
a) A function $\phi :{\cal K}(V) \str \CC$ is called a valuation
if for any $K_1, \, K_2 \in {\cal K}(V)$ such that their union is
also convex one has
$$\phi(K_1 \cup K_2)= \phi(K_1) +\phi(K_2) -\phi(K_1 \cap K_2).$$

b) A valuation $\phi$ is called continuous if it is continuous
with respect to the Hausdorff metric on ${\cal K}(V)$.
\end{definition}

For the classical theory of valuations we refer to the surveys
\cite{mcmullen-schneider} and \cite{mcmullen-survey}. Let us
remind the definition of a {\itshape polynomial} valuation
introduced by Khovanskii and Pukhlikov
\cite{khovanskii-pukhlikov1}, \cite{khovanskii-pukhlikov2}.
\begin{definition}
A valuation $\phi$ is called polynomial of degree at most $d$ if
for every $K\in {\cal K}(V)$ the function $x\mapsto \phi (K+x)$ is
a polynomial on $V$ of degree at most $d$.
\end{definition}
Note that valuations polynomial of degree 0 are called {\itshape
translation invariant} valuations. Polynomial valuations have many
nice combinatorial-algebraic properties
(\cite{khovanskii-pukhlikov1}, \cite{khovanskii-pukhlikov2}).

\begin{example}
(1) The Euler characteristic $\chi$ is a continuous translation
invariant valuation (remind that $\chi(K)=1$ for any convex
compact set $K$).

(2) Let $\mu$ be a measure on $V$ with a polynomial density with
respect to a Lebesgue measure. Fix $A\in {\cal K}(V)$. Then
$$\phi(K):=\mu (K+A)$$ is a continuous polynomial
valuation (here $K+A:=\{k+a|k\in K,\, a\in A\}$).
\end{example}

\begin{remark}
In \cite{alesker-icm} we have announced a result that polynomial
continuous valuations are dense in the space of all continuous
valuations. However during the preparation of this paper we found
a gap in our original argument. So we do not know if this fact is
true.
\end{remark}

Let us remind a basic definition from representation theory. Let
$\rho$ be a continuous representation of a Lie group $G$ in a
Fr\'echet space $F$. A vector $\xi \in F$ is called $G$-smooth if
the map $g\mapsto \rho(g)\xi$ is an infinitely differentiable map
from $G$ to $F$. It is well known (e.g. \cite{wallach}, Section
1.6) that the subset $F^{sm}$ of smooth vectors is a $G$-invariant
linear subspace dense in $F$. Moreover it has a natural topology
of a Fr\'echet space (which is stronger than that induced from
$F$), and the representation of $G$ in $F^{sm}$ is continuous.
Moreover all vectors in $F^{sm}$ are $G$-smooth.

We will denote by $GL(V)$ the group of all linear transformations
of $V$, and by $Aff(V)$ the group of all affine transformations of
$V$.

 We will especially be interested in polynomial
valuations which are $GL(V)$-smooth. The space of $GL(V)$-smooth
valuations polynomial of degree at most $d$ will be denoted by
$PVal_d^{sm}(V)$. This is a Fr\'echet space. Let $PVal^{sm}(V)$
denote the inductive limit of the Fr\'echet spaces
$PVal_d^{sm}(V)$ (with the topology of the inductive limit).

In Section 1 we define a canonical structure of commutative
associative algebra with unit on $PVal^{sm}(V)$ where the unit is
the Euler characteristic. Let us give the idea of this
construction. Let us denote by ${\cal G}'(V)$ the linear space  of
valuations on $\ck(V)$ which are finite linear combinations of
valuations from Example 0.3 (2). It turns out that  ${\cal
G}'(V)\cap PVal^{sm}(V)$ is dense in $PVal^{sm}(V)$ (see the proof
of Lemma 1.1).
 Let $W$ be another linear real vector space. Let
us define the exterior product $\phi \boxtimes  \psi \in {\cal G}'
(V\times W)$ of two valuations $\phi \in{\cal G}'(V), \, \psi\in
{\cal G}'(W)$. Let $\phi (K)=\sum_i \mu_i(K+A_i),\,
\psi(L)=\sum_j\nu_j(L+A_j)$. Define
$$(\phi\boxtimes  \psi) (M):=\sum_{i,j}(\mu_i \boxtimes \nu
_j)(M+(A_i\times  B_j)),$$ where $\mu_i \boxtimes \nu _j$ denotes
the usual product measure.

Now let us define a product on ${\cal G}'(V)$. Let $\Delta:
V\hookrightarrow V\times V$ denote the diagonal imbedding. For
$\phi, \, \psi \in {\cal G}'(V)$ let
$$ \phi \cdot \psi :=\Delta^* (\phi \boxtimes  \psi),$$
where $\Delta ^*$ denotes the restriction of a valuation on
$V\times V$ to the diagonal. By Proposition 1.4 this product is
associative and commutative, the Euler characteristic is the
identity element in this algebra, and this product extends
(uniquely) to a continuous product on $PVal^{sm}(V)$ (Proposition
1.10).

We have the following natural filtration on $PVal^{sm}(V)$:
$$\gamma_i=\{\phi\in PVal^{sm}(V)| \, \phi (K)=0 \, \mbox{ for all } K \mbox{
s. t. } \dim K< i\}.$$ Then $PVal^{sm} (V)=\gamma_0\supset
\gamma_1\supset \dots \supset \gamma_n \supset \gamma_{n+1}=0$.
However this filtration is not compatible with the multiplicative
structure. Nevertheless there exists another filtration $W_i$
which can be characterized as follows.
\begin{theorem}
There exists a unique filtration $W_i$ on $PVal^{sm}(V)$ such that

(1) $\{W_i\}$ is compatible with the multiplicative structure,
i.e. $W_i \cdot W_j \subset W_{i+j}$;

(2) $\gamma_{i+1}\subset W_i\subset \gamma_i \mbox{ for all } i;$

(3) $W_0=\gamma_0,\,W_1 =\gamma_1$;

(4) $W_i$ is a closed subspace of $PVal^{sm}(V)$;

(5) $W_i$ is $Aff (V)$-invariant.

\end{theorem}
This theorem is a reformulation of Theorem 3.8 of Section 3. The
explicit construction of this filtration is given in Section 3.
One of the ways to describe this filtration explicitly is as
follows (see Proposition 3.4):
$$\phi \in W_i \mbox{ if and only if } \lim _{r\str +0} r^{-i+1}\phi(rK+x)=0
\, \mbox{ for all } K\in \ck(V),\, x\in V.$$

We will denote the space of all translation invariant continuous
valuations on $V$ by $Val(V)$. Remind that a valuation $\phi$ is
called $i$-homogeneous if $\phi(\lam K)=\lam^i \phi(K)$ for all
$\lam\geq 0, \, K\in \ck(V)$. We will denote by $Val_i(V)$ the
subspace of $i$-homogeneous valuations. By a result of McMullen
\cite{mcmullen-euler} one has a decomposition:
\begin{equation} Val(V)=\oplus_{i=0}^n Val_i(V).
\end{equation}
Also one has
\begin{equation} Val_0(V)=\CC \cdot \chi,\, Val_n(V)=\CC\cdot
vol\end{equation} where the first equality is trivial, and the
second one is due to Hadwiger \cite{hadwiger-book}. We have also a
further decomposition of these spaces with respect to parity.
Namely we say that a valuation $\phi$ is {\itshape even} if
$\phi(-K)=\phi(K), \, \forall K\in \ck(V)$. Similarly $\phi$ is
called {\itshape odd} if $\phi(-K)=-\phi(K), \, \forall K\in
\ck(V)$. The subspace of even translation invariant valuations
will be denoted by $Val^0(V)$, and the subspace of odd translation
invariant valuations will be denoted by $Val^1(V)$. Similarly
$Val_i^0(V)$ and $Val_i^1(V)$ will denote their subspaces of
$i$-homogeneous valuations. We obviously have
\begin{equation} Val_i(V)=Val_i^0(V)\oplus Val_i^1(V).
\end{equation}

It turns out that one can easily describe the associated graded
algebra $gr_W(PVal^{sm}(V))$ in terms of the algebra of
translation invariant smooth valuations. In Section 3 we prove the
following result (Theorem 3.9).
\begin{theorem}
There exists a canonical isomorphism of graded algebras
$$gr_W (PVal^{sm}(V))\simeq Val^{sm}(V) \otimes \CC[V]$$
where the $i$th graded term in the right hand side is equal to
$Val_i^{sm}(V)\otimes \CC[V]$ where $\CC[V]$ denotes the algebra
of polynomial functions on $V$.
\end{theorem}

In the proof of this theorem we construct this isomorphism
explicitly. Now let us discuss in more detail the case of
translation invariant valuations. One of the main results on
translation invariant valuations we will use in this paper is as
follows (this was proved in \cite{alesker-gafa}).
\begin{theorem}[Irreducibility Theorem]
The natural representations of the group $GL(V)$ in $Val_i^0(V)$
and $Val_i^1(V)$ are irreducible.
\end{theorem}

Below we will denote by $Val^{sm}(V),\, Val_i^{sm}(V),\,
(Val_i^0(V))^{sm},\, (Val_i^1(V))^{sm}$ the subspaces of
$GL(V)$-smooth vectors of the corresponding spaces. The space
$Val^{sm}(V)$ is a subalgebra of $PVal^{sm}(V)$. Moreover the
degree of homogeneity and parity are (obviously) compatible with
the multiplication and define the structure of a bigraded algebra
on $Val^{sm}(V)$ (where the first grading is by $\ZZ$ and the
second grading is by $\ZZ/2\ZZ$). A non-trivial property of the
multiplication is the following version of the Poincar\'e duality
(remind that $Val_n(V)$ is one dimensional) where for a Fr\'echet
space $X$ we denote by $X^*$ its topological dual.
\begin{theorem}
The maps $$(Val_i^0(V))^{sm}\to ((Val_{n-i}^{0}(V))^*)^{sm}\otimes
Val_n(V), \, \mbox{ and }$$ $$ (Val_i^1(V))^{sm}\to
((Val_{n-i}^{1}(V))^*)^{sm}\otimes Val_n(V)$$ induced by the
multiplication $(Val_i^0(V))^{sm}\otimes (Val_{n-i}^0(V))^{sm}\to
Val_n(V)$ and  $(Val_i^1(V))^{sm}\otimes (Val_{n-i}^1(V))^{sm}\to
Val_n(V)$, are isomorphisms.
\end{theorem}
This result is proved in Section 2 (Theorem 2.1). The proof is
based on the Irreducibility Theorem 0.7. In Section 2 we also
compute explicitly some examples of products of valuations
(Proposition 2.2). Note that the isomorphisms in Theorem 0.8
commute with the natural action on the group $GL(V)$. An attempt
to understand these isomorphisms from purely representation
theoretical point of view was done in \cite{alesker-dps}.

If we forget the grading by parity then we obtain a version of the
Poincar\'e duality on the graded algebra
$Val^{sm}(V)=\oplus_{i=0}^n Val_i^{sm}(V).$

We also describe the algebra structure on the space of isometry
invariant continuous valuation on an $n$-dimensional Euclidean
space (as a vector space it is described by the classical Hadwiger
characterization theorem \cite{hadwiger-book}). This graded
algebra is isomorphic to the graded algebra of truncated
polynomials $\CC[x]/(x^{n+1})$ (see Theorem 2.6 for more details).

Now let us discuss some applications of the above results to the
spaces of  valuations invariant under a group. Let $G$ be a
compact subgroup of $GL(V)$. Let us denote by $Val^G(V)$ the space
of $G$-invariant translation invariant continuous valuations.
Again we have a decomposition
$$Val^G(V)=\oplus_{i=0}^n Val_i^G(V).$$
Let us write $h_i:=\dim Val_i^G(V).$
\begin{theorem}
Assume that $G$ is a compact subgroup of $GL(V)$ acting
transitively on the projective space $\PP(V)$.

(i)Then $Val^G(V)$ is finite dimensional.

(ii) $Val^G(V)\subset Val^{sm}(V).$

(iii) $Val^G(V)$ is a finite dimensional graded subalgebra of
$Val^{sm}(V)$ satisfying Poincar\'e duality:
$$Val^G_i(V)\otimes Val_{n-i}^G(V)\str Val_n^G(V)=\CC\cdot vol$$
is a perfect pairing. In particular $h_i=h_{n-i}$.

(iv) $Val_1^G(V)$ is spanned by the intrinsic volume $V_1$, and
$Val^G_{n-1}(V)$ is spanned by the intrinsic volume $V_{n-1}$
where the intrinsic volumes are taken with respect to a
$G$-invariant Euclidean metric. Thus $h_1=h_{n-1}=1$.

(v) Assume in addition that $-Id\in G$. Then the numbers $h_i$
satisfy the Lefschetz inequalities:
$$h_i\leq h_{i+1} \mbox{ for } i<n/2.$$
\end{theorem}
For the definition of the intrinsic volumes we refer to
\cite{schneider-book}. Note that part (i) of this theorem was
proved in \cite{alesker-adv} (Theorem 8.1), and part (ii) in
\cite{alesker-univa} (see the proof of Corollary 1.1.3). Part
(iii) is a direct consequence of Theorem 0.8. Note that the only
proof we know  of the equality $h_i=h_{n-i}$ for a general compact
group $G$ is based on the existence of the multiplicative
structure on the much larger (infinite dimensional) space
$Val^{sm}(V)$ and a version of Poincar\'e duality for this larger
algebra. Note also that Poincar\'e duality for it is based on the
Irreducibility Theorem 0.7. Part (iv) will be proved in Section 2
of this paper. Part (v) of Theorem 0.9 was proved in
\cite{alesker-univa} and it is based on a version of the hard
Lefschetz theorem for even valuations. We expect that this result
(i.e. hard Lefschetz theorem) should be true also for odd
valuations; in that case the condition $-Id\in G$ in part (v) of
Theorem 0.9 could be omitted. However we do not know this.

The paper is organized as follows. In Section 1 we define the
multiplicative structure on polynomial valuations and study its
properties. In Section 2 we study the properties of the subalgebra
of translation invariant smooth valuations. In Section 3 we
introduce and study filtrations on polynomial valuations. In
Section 4 we have some further remarks and discuss some examples.

 {\bf Acknowledgements.} We would like to thank J. Bernstein for very
useful discussions.

\section{The product on valuations.}
Let us agree on a notation. In the rest of the paper we will
denote by $Pol_d(V)$ the space of homogeneous polynomials of
degree $d$. Let us denote by ${\cal G}(V)$ the linear space of
valuations on $V$ which are finite linear combinations of
valuations of the form $K\mapsto \mu(i(K)+A)$ where
$i:V\hookrightarrow V'$ is a linear imbedding of $V$ into a larger
linear space $V'$, and $A\in \ck (V')$ is a fixed convex compact
set, and $\mu $ is a smooth measure on $V'$.

\begin{lemma}
$\cg(V)\cap PVal^{sm}(V)$ is dense in $PVal^{sm}(V)$.
\end{lemma}
{\bf Proof.} Let $\phi $ be a valuation of the above form, i.e.
$\phi(K)=\mu (i(K)+A)$ where $i:V\hookrightarrow V'$, $\mu$ is a
smooth measure on $V'$, and $A\in \ck(V')$. We have to prove that
for any $d$ the space $\cg(V)\cap PVal^{sm}_d(V)$ is dense in
$PVal_d^{sm}(V)$. Let us prove it by induction in $d$. For $d=0$ this
is precisely McMullen's conjecture which was proved in
\cite{alesker-gafa} (as an  easy consequence of the Irreducibility
Theorem 0.7). (Recall that McMullen's conjecture from
\cite{mcmullen-conj} says that the valuations of the form
$K\mapsto vol(K+A)$ are dense in the space $Val(V)$.) Let us
assume that $d>0$. Then we have a continuous map $PVal^{sm}_d(V)
\str Val^{sm}\otimes Pol_d(V)$ which is defined as follows. For
any $\phi \in PVal^{sm}_d(V)$ we have
$$\phi(K+x)=P_d(K)(x)+\mbox{ lower order terms } $$
where $P_d\in (Val(V)\otimes Pol_d(V))^{sm}$. However by Lemma 1.5
below $(Val(V)\otimes Pol_d(V))^{sm}=Val^{sm}(V)\otimes Pol_d(V)$
So the map $\phi \mapsto P_d$ is the desired map. By the
assumption of induction $\cg(V)\cap PVal_{d-1}^{sm}(V)$ is dense in
$PVal_{d-1}^{sm}(V)$. It follows from the case $d=0$ that the image of
$\cg(V)\cap PVal_{d}^{sm}(V)$ in $Val^{sm}(V)\otimes Pol_d(V)$ has a
dense image. Hence the lemma follows. \qed

We need the following technical lemma.
\begin{lemma}
Let $\phi_1,\dots, \phi_s \in \cg(V)$. Then there exists a linear
space $V'$, a linear imbedding $f:V\hookrightarrow V'$, smooth
measures $\mu_i$ on $V'$, and $A_1,\dots, A_s \in \ck(V')$ such
that $\phi_i(K)=\mu_i(f(K)+A_i),\, i=1,\dots ,s$.
\end{lemma}
{\bf Proof.} By definition $\phi_i$ has the form
$$\phi_i(K)=\nu_i(g_i(K)+B_i)$$
where $g_i:V\hookrightarrow V_i',\, B_i\in \ck(V_i')$, and $\nu_i$
are smooth measures on $V_i'$. Let $V'$ be the inductive limit of
the system $\{V\overset{g_i}{\hookrightarrow} V_i'\}_{i=1}^s$. It
would be convenient to have an explicit construction for it. Let
$T_i$ be a complement of $g_i(V)$ in $V_i'$. Then $V'$ is
isomorphic to $V\oplus T_1\oplus \dots \oplus T_s$. With this
decomposition $f$ is the identity imbedding. Let us construct
$A_i,\, \nu_i$ say for $i=1$. Note that $V_1'$ is isomorphic to
$V\oplus T_1$. Let us denote also by $B_1$ the image of $B_1$ in
$V'$. Let $S:=V_2\oplus \dots \oplus V_s$. Let us fix a Lebesgue
measure $dvol_S$ on $S$. Fix $Q\in \ck(S)$ any set of volume 1.
Set $A_1:=B_1\times Q$. Let $\mu_1:=\nu_1 \boxtimes dvol_S$. These
choices satisfy the lemma. \qed

 Let $W$ be another linear real vector space. Let
us define the exterior product $\phi \boxtimes  \psi \in {\cal G}
(V\times W)$ of two valuations $\phi \in{\cal G}(V), \, \psi\in
{\cal G}(W)$. By Lemma 1.2 we may assume that these valuations
have the form $\phi (K)=\sum_i \mu_i(f(K)+A_i),\,
\psi(L)=\sum_j\nu_j(g(L)+B_j)$ where $f:V\hookrightarrow V',\,
g:W\hookrightarrow W'$ are imbeddings, $A_i\in \ck(V'),\, B_j\in
\ck(W')$. Define
$$(\phi\boxtimes  \psi) (M):=\sum_{i,j}(\mu_i \boxtimes \nu
_j)((f\times g)(M)+(A_i\times B_j)),$$ where $\mu_i \boxtimes \nu
_j$ denotes the usual product measure on $V'\times W'$.

Let us denote by $CVal(V)$ the closure of $\cg(V)$ in the
Fr\'echet space of all continuous valuations on $V$ with the
topology of uniform convergence on compact subsets on $\ck(V)$.
\begin{proposition}
(i) For $\phi \in{\cal G}(V), \, \psi\in {\cal G}(W)$ their
exterior product $\phi\boxtimes  \psi \in {\cal G}(V\times W)$ is
well defined.

(ii) The exterior product is bilinear   with respect to each
argument.

(iii) Fix $\phi \in \cg (V)$. Then the map $\cg (W)\str \cg
(V\times W)$ given by $\psi\mapsto \phi \boxtimes \psi$ extends
(uniquely) by continuity to a map $CVal(V)\str CVal(V\times W)$.

(iv)
$$(\phi \boxtimes \psi)\boxtimes \eta= \phi \boxtimes (\psi\boxtimes
\eta).$$

(v) Let $f:V\hookrightarrow V_1,\, g:W\hookrightarrow W_1$ be two
imbeddings. Let $\phi \in \cg(V_1),\, \psi\in \cg(W_1)$. Then
$$(f\times g)^*(\phi\boxtimes \psi)=f^*\phi \boxtimes g^*\psi.$$

\end{proposition}
{\bf Proof.} Note that the parts (ii), (iv), and (v) are obvious.
Now let us fix a valuation $\phi \in \cg (V)$ of the form $\phi
(K)=\mu(K+A),\, A\in V'$. From the definition of the exterior
product and the Fubini theorem one easily gets the following
formula:
$$(\phi\boxtimes \psi) (M)= \int _{x\in V'}\psi \left(((f\times g)(M)+(f(A)\times \{0\}))\cap
(\{x\} \times W)\right)d\mu (x).$$ Thus the right hand side in the
above formula does not depend on the particular form of
presentation of $\psi$. Now the parts (i) and (iii) follow easily.
\qed

Now let us define a product on ${\cal G}(V)$. Let $\Delta:
V\hookrightarrow V\times V$ denote the diagonal imbedding. For
$\phi, \, \psi \in {\cal G}(V)$ let
$$ \phi \cdot \psi :=\Delta^* (\phi \boxtimes  \psi),$$
where $\Delta ^*$ denotes the restriction of a valuation on
$V\times V$ to the diagonal.
\begin{proposition}
Equipped with the above defined multiplication, $\cg(V)$ becomes
an associative commutative unital algebra where the unit is the
Euler characteristic $\chi$.
\end{proposition}
{\bf Proof.}  The associativity follows from  Proposition 1.3
(iv). The commutativity is obvious. Let us prove that the Euler
characteristic $\chi$ is the unit. Let a valuation $\psi$ have the
form $\psi(K)=\mu(K+A)$ where $A$ is a fixed set from $\ck(V)$.
Let $\Delta: V\str V\times V$ be the diagonal imbedding. Then by
the definition of the product and the Fubini theorem we have
$$(\chi \cdot \psi) (K)= \int _{x\in V}\chi \left((\Delta K+(A\times \{0\}))\cap
(\{x\} \times V)\right)d\mu (x)=\mu(K+A)=\psi(K).$$ \qed

Next let us prove the following technical lemma which is well
known.
\begin{lemma}
Let ${\cal F}$ be a Fr\'echet $G$-module. Let $S$ be a finite
dimensional $G$-module. Then $({\cal F}\otimes S)^{sm}={\cal
F}^{sm}\otimes S$.
\end{lemma}
{\bf Proof.} It is well known and easy to see that for any finite
dimensional $G$-module $S$ one has $S^{sm}=S$. Hence
$\cf^{sm}\otimes S\subset (\cf\otimes S)^{sm}$. To prove the
opposite inclusion let us fix $\xi\in(\cf\otimes S)^{sm}$. Let us
also fix a basis $s_1,\dots,s_k$ of $S$. Then
$\xi=\sum_{i=1}^kf_i\otimes s_i$ with $f_i\in \cf$. We have to
show that $f_i\in \cf^{sm}$. Let $\{s_i^*\}$ be the dual basis in
$S^*$. One has the canonical map $t:\cf\otimes S\otimes S^*\str
\cf$ given by $t(f\otimes s\otimes s^*)=s^*(s)\cdot f$. Note that
$t(f\otimes s_i\otimes s_i^*)=f$. Moreover $t$ commutes with the
action of $G$. Hence $t((\cf\otimes S\otimes S^*)^{sm})\subset
\cf^{sm}$. We also have $(\cf\otimes S)^{sm}\otimes S^*\subset
(\cf\otimes S\otimes S^*)^{sm}$. Hence $t((\cf\otimes
S)^{sm}\otimes S^*)\subset \cf^{sm}$. But $f_i=t(\xi\otimes
s_i^*)$. Hence $f_i\in \cf^{sm}$.
 \qed

Let us also recall two well known results we are going to use,
namely the Casselman-Wallach theorem and the L. Schwartz kernel
theorem.

\begin{theorem}[Casselman-Wallach,
\cite{casselman}]\label{casselman-wallach}
 Let $G$ be a real
reductive Lie group. Let $(\rho,G,\cf)$ and $(\pi,G,\cg)$ be
continuous representations of $G$ of moderate growth in Fr\'echet spaces $\cf$ and
$\cg$. Assume in addition that $\cg$ is an admissible $G$-module
of finite length, and
$$\cf^{sm}=\cf,\, \cg^{sm}=\cg.$$
Then any continuous morphism of $G$-modules $f:\cf\to \cg$ has
closed image.
\end{theorem}

\begin{theorem}[L. Schwartz kernel theorem,
\cite{gelfand-vilenkin}]\label{schwartz}
 Let $X_1$and $X_2$ be compact smooth
manifolds. Let $\ce_1$ and $\ce_2$ be smooth finite dimensional
vector bundles over $X_1$ and $X_2$ respectively. Let $\cg$ be a
Fr\'echet space. Let
$$B:C^\infty(X_1,\ce_1)\times C^\infty(X_2,\ce_2)\to \cg$$
be a continuous bilinear map. Then there exists unique continuous
linear operator
$$b:C^\infty(X_1\times X_2,\ce_1\boxtimes \ce_2)\to \cg$$
such that $b(f_1\otimes f_2)=B(f_1,f_2)$ for any $f_i\in
C^\infty(X_i,\ce_i),\, i=1,2$.
\end{theorem}

 Let $\pvd$ denote the
space of continuous valuations on $V$ which are polynomial of
degree at most $d$. Let $\ond$ denote the (finite dimensional)
space of $n$-densities on $V$ with polynomial coefficients of
degree at most $d$ (clearly $\ond$ is canonically isomorphic to
$(\oplus _{i=0}^{d}Sym^{i}V^*) \otimes |\wedge ^{n}V^*|$ where
$|\wedge ^nV^*|$ denotes the space of Lebesgue measures on $V$).

Let us denote by $\pl$ the manifold of oriented lines passing
through the origin in $V^*$. Let $L$ denote the line bundle over
$\pl$ whose fiber over an oriented line $l$ consists of linear
functionals on $l$.

We are going to construct  a natural linear map
$$\Theta_{k,d}: \ond \otimes C^{\infty} ((\pl)^k, L^{\boxtimes
k})\str \pvd$$ which commutes with the natural action of the group
$GL(V)$ on both spaces and induces an epimorphism on the subspaces
of smooth vectors.

The construction is as follows. Let $\mu\in \ond,\, A_1,\dots
,A_k\in \ck(V)$. Then $\int_{\sum_{j=1}^k \lam_jA_j} \mu$ is a
polynomial in $\lam_j\geq 0$ of degree at most $n+d$. This can be
easily seen directly, but it was also proved in general for
polynomial valuations by Khovanskii and Pukhlikov
\cite{khovanskii-pukhlikov1}. Also it easily follows that the
coefficients of this polynomial depend continuously on $(A_1,\dots
,A_k)\in \ck(V)^k$ with respect to the Hausdorff metric. Hence we
can define a continuous map $\Theta_{k,d}':\ond \times\ck(V)^k\to
PVal_d(V)$ given by
$$(\Theta_{k,d}'(\mu;A_1,\dots,A_k))(K):=\frac{\pr^k}{\pr
\lam_1\dots\pr\lam_k}\bigg |_{\lam_j =0}\int_{K+\sum_{j=1}^k
\lam_jA_j}\mu.$$

It is clear that $\Theta_{k,d}'$ is Minkowski additive with
respect to each $A_j$. Namely, say for $j=1, \, a,b\geq 0$, one
has
$$\Theta_{k,d}'(\mu;a A_1'+b
A_1'',A_2,\dots,A_k)=a\Theta_{k,d}'(\mu;A_1',A_2,\dots,A_k)+
b\Theta_{k,d}'(\mu;A_1'',A_2,\dots,A_k).$$

Remind that for any $A\in \ck(V)$ one defines the supporting
functional $h_A(y):=\sup_{x\in A}(y,x)$ for any $y\in V^*$. Thus
$h_A\in C(\PP_+(V^*),L)$. Moreover it is well known (and easy to
see) that $A_N\str A$ in the Hausdorff metric if and only if
$h_{A_N}\str h_A$ in $C(\PP_+(V^*),L)$. Also any section $F\in
C^2(\PP_+(V^*),L)$ can be presented as a difference $F=G-H$ where
$G,\,H\in C^2(\PP_+(V^*),L)$ are supporting functionals of some
convex compact sets and $\max\{||G||_2,||H||_2\}\leq c ||F||_2$
where $c$ is a constant. (Indeed one can choose $G=F+R\cdot h_D,\,
H=R\cdot h_D$ where $D$ is the unit Euclidean ball, and $R$ is a
large enough constant depending on $||F||_2$.) Hence we can
uniquely extend $\Theta_{s,d}'$ to a multilinear continuous map
(which we will denote by the same letter):
$$\Theta_{k,d}':\ond \times (C^2(\PP_+(V^*),L))^k\to PVal_d(V).$$
By the L. Schwartz kernel theorem (Theorem \ref{schwartz}) it
follows that this map gives rise to a continuous linear map
$$\Theta_{k,d}:\ond\otimes C^{\infty}(\PP_+(V^*)^k, L^{\boxtimes
k})\str PVal_d^{sm}(V)$$ which we wanted to construct.





We will study this map $\Theta_{k,d}$. Note that it depends on $k$
and $d$ which will be fixed from now on. Let us denote by
$\Theta_d$ the sum of the maps $\bigoplus
_{k=0}^{n-1}\Theta_{k,d}$. Thus $\Theta_d:\ond \otimes(\bigoplus
_{k=0}^{n-1}C^{\infty}(\pl^k, L^{\boxtimes k})) \str
PVal^{sm}_d(V).$
\begin{lemma}
The map $\Theta_d$  has a dense image.
\end{lemma}
\begin{corollary}
The map $\Theta_d$ is an epimorphism.
\end{corollary}
This corollary follows immediately from Lemma 1.8 and the
Casselman-Wallach theorem (Theorem \ref{casselman-wallach}). So
let us prove Lemma 1.8.

{\bf Proof} of Lemma 1.8. We will prove it by induction in $d$.
For $d=0$ this is just McMullen's conjecture proved in
\cite{alesker-gafa}. Assume that the statement of Lemma 1.8 holds
for $d-1$. It is sufficient to show that the map
$Sym^d(V^*)\otimes |\wedge^{n} V^*| \otimes(\bigoplus
_{k=0}^{n-1}C^{\infty}(\pl^k, L^{\boxtimes k})) \to
(\pvd/PVal_{d-1}(V))^{sm}$ is onto. However the last space is
equal to $(Sym^d(V^*)\otimes Val(V))^{sm}=Sym^d(V^*)\otimes
Val^{sm}(V)$ (using Lemma 1.5). Thus we have reduced the claim
again to the case $d=0$. \qed

\begin{proposition}
The multiplication $$(PVal_i(V))^{sm}\otimes (PVal_j(V))^{sm}\str
(PVal_{i+j}(V))^{sm}$$ is continuous.
\end{proposition}
{\bf Proof.} Since this map commutes with the action of $GL(V)$ it
is sufficient to prove that the map
$$(PVal_i(V))^{sm}\otimes (PVal_j(V))^{sm}\str
PVal_{i+j}(V)$$ is continuous. In fact we will prove  a bit more.
We will prove that the exterior product of valuations is a
continuous bilinear map $$(PVal_i(V))^{sm}\times
(PVal_j(W))^{sm}\str PVal_{i+j}(V\times W).$$ For the simplicity
of notation we will assume that $W=V$.

Let us fix a large natural number $N\geq 2$. Let $\phi \in
(PVal_i(V))^{sm},\, \psi \in (PVal_j(V))^{sm}$ be such that $\phi=
\sum _p \omega _p\otimes l_p,\, \psi =\sum _q \omega'_q\otimes
l'_q$ where $\omega _p\in \Omega^{n}_{i},\, \omega'_q\in
\Omega^{n}_{j},\, l_p,\, l'_q \in \bigoplus
_{k=0}^{n-1}C^{\infty}(\pl^k, L^{\boxtimes k})$ and
$$max\{\sum_p||\omega _p||,\, \sum_q||\omega'_q||,\,
\sum_p||l_p||_{2N},\, \sum_q||l'_q||_{2N}\}\leq 1$$ where $||\cdot
||$ is a norm on $\Omega^n_i$, and $||\cdot||_N$ is a norm on
$\bigoplus _{k=0}^{n-1}C^{2N}(\pl^k, L^{\boxtimes k}).$

 Then it is
easy to see that for a constant $C$ (depending on $n$ and $k$
only) such that one can write $$l_p=\sum_{s=1}^{\infty}
\gamma_{s,p,1}\otimes \dots \otimes\gamma_{s,p,k},\,
l'_q=\sum_{s=1}^{\infty}\gamma'_{s,q,1}\otimes \dots
\otimes\gamma'_{s,q,k}$$ where $\gamma_{s,p,i},\,\gamma
_{s,q,i}\in C^{\infty}(\pl, L)$ and
$$max\{\sum_s \sum_{k=0}^{n-1}\prod_{t=1}^{k}||\gamma_{s,p,t}||_N,\,
\sum_s\sum_{k=0}^{n-1}\prod_{t=1}^k||\gamma'_{s,q,t}||_N \}<C$$
where $C$ is a constant which might be different from the previous
$C$. Every $\gamma \in C^{\infty}(\pl, L)$ can be written as
$\gamma= \gamma_1 -\gamma _2$ such that $\gamma_1,\, \gamma _2$
are supporting functionals of convex sets and
$max\{||\gamma_1||_N, ||\gamma_2||_N\}\leq C||\gamma||_N$ (using
the same argument as in the construction of the maps
$\Theta_{k,d}$). Hence we may assume that all $\gamma_{s,p,t},\,
\gamma'_{s,q,t}$ are supporting functionals of convex sets
$A_{s,p,t},\, A'_{s,q,t}$ respectively. Let
$\Omega_{p,q}:=\omega_p \wedge \omega'_q$; it is a top-degree form
on $V\times V$. Then by the definition of the exterior product of
valuations
$$(\phi \boxtimes \psi)(K) =$$
$$\sum_{s=1}^\infty \sum_{k,k'=0}^{n-1}\sum_{p,q} \frac
{\partial^{k}}{\partial \lambda_1 \dots \partial \lambda_{k}}
\frac {\partial^{k'}}{\partial \mu_1 \dots
\partial \mu_{k'}}\bigg |_0\int _{K+\sum_{m=1}^k\lambda_m (A_{s,p,m}\times 0)+\sum_{m'=1}^{k'}
\mu_{m'}(0\times A'_{s,p,m'})} \Omega_{p,q} .$$ Note that
$$\frac
{\partial^{k}}{\partial \lambda_1 \dots \partial \lambda_{k}}
\frac {\partial^{k'}}{\partial \mu_1 \dots
\partial \mu_{k'}}\bigg|_0\int _{K+\sum_{m=1}^k\lambda_m (A_{s,p,m}\times 0)+\sum_{m'=1}^{k'}
\mu_{m'}(0\times A'_{s,q,m'})} \Omega_{p,q}= $$
$$\prod_m ||\gamma_{s,p,m}||_0 \prod_{m'}||\gamma_{s,q,m'}||_0\cdot $$
$$\frac {\partial^{k}}{\partial \lambda_1 \dots \partial
\lambda_{k}} \frac {\partial^{k'}}{\partial \mu_1 \dots
\partial \mu_{k'}}\bigg |_0\int _{K+\sum_{m=1}^k\lambda_m (\frac{A_{s,p,m}}{||\gamma_{s,p,m}||_0}\times 0)+
\sum_{m'=1}^{k'} \mu_{m'}
(0\times\frac{A'_{s,q,m'}}{||\gamma_{s,q,m'}||_0})}
\Omega_{p,q}.$$

The last integral is a polynomial in $\lambda_m,\, \mu_{m'}$ of
degree at most $i+j+2n$. Hence in order to get an estimate on its
derivatives it is sufficient to estimate the polynomial itself for
all $\lam _m,\,\mu_{m'}$ of absolute value at most 1. If the set
$K$ is uniformly bounded then the whole Minkowski sum of sets
under the integral is uniformly bounded. Hence all the integrals
are uniformly bounded, and we get
$$|(\phi \boxtimes \psi)(K)|\leq C\sum_{s=1}^\infty \sum_{k,k'=0}^{n-1}\sum_{p,q}
\left( \prod_m ||\gamma_{s,p,m}||_0
\prod_{m'}||\gamma_{s,q,m'}||_0\right).$$ The last sum is
estimated by a constant. \qed

\section{Translation invariant valuations.} In this section we
discuss the subalgebra $(Val(V))^{sm}$ of smooth translation
invariant valuations. It has a structure of a Frobenius algebra,
namely it satisfies a version of the Poincar\'e duality with
respect to the natural grading.

Remind that by a result of Hadwiger \cite{hadwiger-book} the space
$Val_n(V)$ is one dimensional and it is spanned by  Lebesgue
measure. First observe that the multiplication
$(Val_i^{0}(V))^{sm}\otimes (Val_{n-i}^{1}(V))^{sm}\str Val_n(V)$
is trivial since the product of even and odd valuations must be
odd, and all valuations of the maximal degree of homogeneity are
even. The main result of this section is
\begin{theorem}
(i) The map $(Val_i^{0}(V))^{sm}\otimes
(Val_{n-i}^{0}(V))^{sm}\str Val_n(V)$ is a perfect pairing. More
precisely the induced map
$$(Val_i^{0}(V))^{sm}\str (Val_{n-i}^{0}(V)^*)^{sm}\otimes
Val_n(V)$$ is an isomorphism.

(ii) For $1\leq i\leq n-1$, the map $(Val_i^{1}(V))^{sm}\otimes
(Val_{n-i}^{1}(V))^{sm}\str Val_n(V)$ is a perfect pairing in the
above sense.
\end{theorem}
Before we prove the theorem we will need a proposition which is of
independent interest.
\begin{proposition}
Let $V$ be an $n$-dimensional Euclidean space. Let
$$\phi(K)=V(K[i],A_1,\dots, A_{n-i}),\,
\psi(K)=V(K[n-i],B_1,\dots,B_i)$$ where $A_p,\, B_q \in \ck (V)$
are fixed. Then
$$(\phi \cdot \psi)(K)={n\choose i}^{-1}V(A_1,\dots,
A_{n-i},-B_1,\dots, -B_{i}) vol(K).$$
\end{proposition}

First we have the following simple identity.
\begin{claim}
$$V(K[i],A_1,\dots, A_{n-i})=(n(n-1)\dots
(i+1))^{-1}\frac{\pr^{n-i}}{\pr \lam_1 \dots \pr \lam_{n-i}}\bigg
|_{0} vol (K+\sum _{j=1}^{n-i} \lam_j A_j).$$
\end{claim}
{\bf Proof} of Proposition 2.2. Using Claim 2.3 we have
$$(\phi \cdot \psi)(K)=(n(n-1)\cdots
(i+1))^{-1}(n(n-1)\cdots (n-i+1))^{-1}\cdot$$
$$\frac{\pr^{n-i}}{\pr \lam_1 \cdots \pr \lam_{n-i}}\bigg |_{0}
\frac{\pr^{i}}{\pr \mu_1 \cdots \pr \mu_{i}}\bigg |_0
vol_{2n}(\Delta (K)+\sum_{j=1}^{n-i} \lam_j (A_j\times 0)
+\sum_{l=1}^{i}\mu_l (0\times B_l))$$ where $\Delta:
V\hookrightarrow V\times V$ is the diagonal imbedding. Again by
Claim 2.3
$$\frac{\pr^{n-i}}{\pr \lam_1 \cdots \pr \lam_{n-i}}\bigg |_{0}
\frac{\pr^{i}}{\pr \mu_1 \cdots \pr \mu_{i}}\bigg |_0
vol_{2n}(\Delta (K)+\sum_{j=1}^{n-i} \lam_j (A_j\times 0)
+\sum_{l=1}^{i}\mu_l (0\times B_l))=$$
$$(2n\cdots (n+1))\cdot
V_{2n}(\Delta(K)[n]; A_1\times 0,\dots, A_{n-i}\times 0; 0\times
B_1, \dots, 0\times B_i).$$ Hence we obtain
$$(\phi \cdot \psi)(K)={2n \choose n}{n\choose i}^{-1}
V_{2n}(\Delta(K)[n]; A_1\times 0,\dots, A_{n-i}\times 0; 0\times
B_1, \dots, 0\times B_i).$$ We will need a lemma.
\begin{lemma}
Let $X=Y\oplus Z$ be an orthogonal decomposition of a Euclidean
space $X$. Let $\dim X=N,\, \dim Y=n$. Let $M\in \ck(Y),\,
A_1,\dots A_{N-n}\in \ck(X).$ Then
$$V_N(M[n]; A_1,\dots A_{N-n})={N\choose n}^{-1}
vol_n(M)V_{N-n}(Pr_ZA_1,\dots ,Pr_ZA_{N-n})$$ where $Pr_Z$ denotes
the orthogonal projection onto $Z$.
\end{lemma}
Let us postpone the proof of this lemma and continue proving
Proposition 2.2. In our situation $X=V\times V,\,
Y=\Delta(V)=\{(x,x)\},\, Z=\{(x,-x)\}.$ Note that
$$Pr_Z((x,0))=\frac{1}{2}(x,-x),\, Pr_Z((0,x))=\frac{1}{2}(-x,x).$$
Let us denote by $\Delta '$ the imbedding $V\hookrightarrow
V\times V$ given by $\Delta'(x)=(x,-x)$. Using Lemma 2.4 we get
\begin{eqnarray*}
(\phi\cdot \psi)(K)&=& {n\choose i}^{-1} vol_n(\Delta
(K))V_n(\frac{\Delta 'A_1}{2},\dots ,\frac{\Delta
'A_{n-i}}{2},\frac{-\Delta 'B_1}{2}, \dots, \frac{-\Delta
'B_i}{2})\\
&=& {n\choose i}^{-1}V(A_1,\dots, A_{n-i},-B_1,\dots, -B_{i})
vol(K).\end{eqnarray*} \qed

{\bf Proof} of Lemma 2.4. We may assume that $A_1=\dots =
A_{N-n}=A$. We have
$$V_N(M[n],A[N-n])=\frac{(N-n)!}{N!} \frac{d^n}{d\eps ^n}\bigg |_0
vol _N(A+\eps M).$$ We have
\begin{eqnarray*}
vol_{N}(A+\eps M) &=&\int _{z\in Z} vol_n((A+\eps M)\cap
(z+Y))dz\\
&=&\int _{z\in Z} vol_n((A\cap(z+Y))+\eps M)dz\\
&=&\int _{z\in Pr_Z A}(\eps ^n vol _n M +O(\eps ^{n-1}))dz\\
&=&\eps^n vol_{N-n}(Pr_Z A) vol_n M +O(\eps^{n-1}).
\end{eqnarray*}
This proves Lemma 2.4. \qed


{\bf Proof} of Theorem 2.1.  The kernel of the
 map $$(Val_i^{0}(V))^{sm}\str
(Val_{n-i}^{0}(V)^*)^{sm}\otimes Val_n(V)$$ is $GL(V)$-invariant
subspace (and similarly in the odd case). Hence by the
Irreducibility Theorem 0.7 it must be either trivial or coincide
with the whole space. We will show that the multiplication in
non-trivial. Assume that it is proved. Then by the Irreducibility
Theorem 0.7 the image the above map is a dense subspace, and hence
it is equal to the whole space by the Casselman-Wallach theorem.
Thus it remains to check in both cases that the product is
non-zero.

First let us check it in the even case. By Proposition 2.2 the
product of the intrinsic volumes $V_i\cdot V_{n-i}\ne 0$.

Now let us consider the odd case. We want to show that the
multiplication $(Val_i^{1}(V))^{sm}\otimes
(Val_{n-i}^{1}(V))^{sm}\str Val_n(V)$ is non-zero. We have
$$(Val^{\eps}_i(V))^{sm}\otimes (Val^{\delta}_j(V))^{sm}\str
(Val^{\eps +\delta}_{i+j}(V))^{sm}$$ where the addition of upper
indexes is understood modulo 2.

\begin{lemma}
The multiplication $$(Val^{\eps}_1(V))^{sm}\otimes
(Val^{\eps}_1(V))^{sm}\str (Val^{0}_{2}(V))^{sm}$$ is non-zero for
$\eps\in \ZZ/2\ZZ$ and hence has a dense image.
\end{lemma}

First let us finish proving Theorem 2.1 assuming this lemma.
Consider the multiplication maps (we will omit for brevity the
upper sign $sm$):
$$Val^1_1(V) \otimes (Val^0_1(V))^{\otimes (i-1)}\str Val_i^1(V);$$
$$Val^1_1 \otimes (Val^0_1)^{\otimes (n-1-i)}\str Val_{n-i}^1.$$
It is sufficient to showing that the composition
$$\left(Val^1_1(V) \otimes (Val^0_1(V))^{\otimes (i-1)}\right)\otimes
\left(Val^1_1(V) \otimes (Val^0_1(V))^{\otimes (n-1-i)}\right)\str
Val_n(V)$$ is non-zero. This is equivalent to show that
$$Val^1_1(V)\otimes Val^1_1(V) \otimes (Val^0_1(V))^{\otimes (n-2)}\str Val_n(V)$$
is non-zero. Note that the image $(Val^0_1(V))^{\otimes(n-2)}\str
Val_{n-2}^0(V) $ is non-zero since the power of the intrinsic
volume $(V_1)^{n-2}$ does not vanish. Hence this image is a dense
subspace in $ Val_{n-2}^0(V)$. Hence it is sufficient to show that
the map $$Val^1_1(V)\otimes Val^1_1(V)\otimes Val_{n-2}^0(V) \str
Val_n(V)$$ is non-zero. By Lemma 2.5 the map $Val^1_1(V)\otimes
Val^1_1(V)\str Val^0_2(V) $ has a dense image. As we have
previously proved the map $Val^0_2(V)\otimes Val_{n-2}^0(V)\str
Val_n(V) $ is non-zero. This immediately implies Theorem 2.1. \qed

{\bf Proof} of Lemma 2.5. First let reduce the statement to the
case $n=2$. Let us fix a 2-dimensional subspace $W\subset V$. It
is easy to see that the restriction map $Val(V)\str Val(W)$ is
non-zero and its image is $GL(W)$-invariant. Hence by the
Irreducibility Theorem the image of this restriction is a dense
subspace. Moreover the multiplication commutes with the
restriction by Proposition 1.3(v). Hence it is sufficient to show
that
$$Val^1_1(W)\otimes Val^1_1(W)\str Val_2(W)$$
is non-zero.  By Proposition 2.2 we have  for
 $A,B\in  \ck(W)$:
$$V(\cdot, A)\cdot V(\cdot, B)=\kappa V(A,-B) \cdot vol(\cdot),$$
where $\kappa$ is a non-zero normalizing constant.

Now let us choose valuations $\phi(\cdot):=V(\cdot, A)-V(\cdot,
-A),\, \psi(\cdot):=V(\cdot, B)-V(\cdot, -B)$ with $A$ and $B$ to
be chosen later. Then we get
$$(\phi \cdot \psi)(K)=2\kappa (V(A,B)-V(A,-B))\cdot vol(K).$$
If we choose $A$ and $B$ so that $V(A,B)\ne V(A,-B)$ then the
above product does not vanish.\qed

For any subgroup $G$ of the group of linear transformations of the
space $V$ let us denote by $Val^G(V)$ the space of translation
invariant continuous valuations invariant with respect to $G$. Let
now $V$ be a Euclidean space of dimension $n$. Also $Val^G_i(V)$
denote the subspace of $Val^G(V)$ of $i$-homogeneous valuations.
Then it immediately follows from McMullen's theorem
\cite{mcmullen-euler} that $Val^G(V)=\oplus_{i=0}^n Val_i^G(V)$.
Let $O(n)$ denote the full orthogonal group, and $SO(n)$ denote
the special orthogonal group. Let $V_i$ denote the $i$th intrinsic
volume (see e.g. \cite{schneider-book}). It was proved by Hadwiger
that for any $i$ $Val_i^{O(n)}=Val_i^{SO(n)}=\CC \cdot V_i$. Now
we describe the algebra structure on the space $Val^{O(n)}(V)$.
\begin{theorem}
There exists an isomorphism of graded algebras $\CC[x]/(x^{n+1})
\simeq Val^{O(n)}(V)$ given by $x\mapsto V_1$.
\end{theorem}
{\bf Proof.} For any $i$ the valuation $(V_1)^i\in Val_i^{O(n)}$.
Hence it must be proportional to $V_i$. By Proposition 2.2 the
constant of proportionality does not vanish. This implies the
result. \qed

Now let us prove part (iv) of Theorem 0.9. Namely let us prove
that if $G$ is a compact subgroup of $GL(V)$ acting transitively
on the projective space $\PP_+(V)$ then $Val_1^G(V)=\CC\cdot V_1$
and $Val_{n-1}^G=\CC\cdot V_{n-1}$. Here we assume that we
consider the intrinsic volumes $V_i$ with respect to a Euclidean
metric on $V$ invariant under $G$ (which is unique up to a
proportionality). It was proved by McMullen \cite{mcmullen-conj}
that every translation invariant continuous $(n-1)$-homogeneous
valuation $\phi$ has the form
$$\phi(K)=\int_{S^{n-1}} f(x) dS_{n-1}(K,x)$$
where $dS_{n-1}(K,\cdot)$ denotes the $(n-1)$-area measure of $K$
(see \cite{schneider-book}), and $f$ is a continuous function on
the sphere $S^{n-1}$. Moreover $f$ can be chosen to be orthogonal
(with respect to the Lebesgue measure on $S^{n-1}$) to any linear
functional on the sphere, and under this restriction $f$ is
defined uniquely by $\phi$. Hence in our situation $f$ can be
chosen $G$-invariant.  If $G$ acts transitively on the projective
space then it acts transitively on the sphere. Hence $f$ must be
constant. Hence $\phi\in \CC \cdot V_{n-1}$.

Now let us consider the case of $ Val_1^G(V)$. The result follows
from the previous case and the Poincar\'e duality. A more
elementary way to see it is as follows. Let us construct a
canonical imbedding $Val_1(V)\hookrightarrow C^{-\infty}(\PP_+(V),
L^*\otimes |\omega|) )$ where the target space is the space of
generalized sections of the line bundle $L^*\otimes |\omega|$ over
$\PP_+(V)$, where $L$ was defined in Section 1, and  $|\omega |$
is the line bundle of densities over $\PP_+(V)$. This imbedding
(in fact in a more general situation) was first considered by
Goodey and Weil \cite{goodey-weil}; now we essentially reproduce
their argument. It is well known that any valuation $\phi\in
Val_1(V)$ is Minkowski additive, i.e. $\phi(\lam A+\mu B)=\lam
\phi(A)+\mu \phi(B)$ for all $A,\, B\in \ck(V)$ and $\lam, \, \mu
\geq 0$. Note that the supporting functional of any convex compact
set is a continuous section of the line bundle $L$ over $\PP_+(V)$
defined in Section 1. However any $C^{\infty}$-section $F$ of $L$
can presented as a difference of two smooth supporting functionals
of convex compact sets, $F=G-H$ such that
$\max\{||G||_2,||H||_2\}\leq ||F||_2$. Hence using the Minkowski
additivity of $\phi$ we can extend it uniquely to a continuous
linear functional on $C^{\infty}(\PP_+(V),L)$. Clearly we get a
continuous imbedding $Val_1(V)\hookrightarrow
C^{-\infty}(\PP_+(V), L^*\otimes |\omega|) ).$ However if the
group $G$ acts transitively on the sphere then the space of
$G$-invariant (generalized) sections is at most one dimensional.
It follows that $Val_1^G(V)=\CC \cdot V_1$. \qed

\section{Filtrations on polynomial valuations.}
Let us define the following filtration on $PVal^{sm}(V)$:
$$\gamma_i=\{\phi\in PVal^{sm}(V)| \, \phi (K)=0 \, \forall K \mbox{
such that } \dim K< i\}.$$ Then $PVal^{sm} (V)=\gamma_0\supset
\gamma_1\supset \dots \supset \gamma_n \supset
\gamma_{n+1}=\{0\}$. Note that $$\gamma_i\cap
Val^{sm}(V)=(Val^1_{i-1}(V))^{sm}\oplus Val_i^{sm}(V)\oplus \dots
\oplus Val_n^{sm}(V).$$

Let us introduce another filtration on $(PVal)^{sm}$. Set
$$W_i(PVal_d^{sm}):=\sum _{k=0}^{n-i}Im(\Theta_{k,d})$$
where the maps $\Theta_{k,d}$ were defined in Section 1. We have
$$PVal_d^{sm}=W_0(PVal_d^{sm})\supset W_1(PVal_d^{sm})\supset
\dots \supset W_n(PVal_d^{sm}).$$ Moreover $W_n(PVal_d^{sm})$
coincides with densities on $V$ polynomial of degree at most $d$.

\begin{lemma} For any $d$
$$W_i(PVal_d^{sm})\cap Val^{sm}(V)=Val_i^{sm}(V)\oplus Val_{i+1}^{sm}(V)\oplus \dots
\oplus Val_n^{sm}(V).$$
\end{lemma}
{\bf Proof.} Obviously we have an inclusion
$$Val_i^{sm}(V)\oplus Val_{i+1}^{sm}(V)\oplus \dots
\oplus Val_n^{sm}(V)\subset W_i(PVal_d^{sm})\cap Val^{sm}(V).$$
Since $W_i(PVal_d^{sm})\subset \gamma_i$ (by Proposition 3.7(i)
below) we have
$$W_i(PVal_d^{sm})\cap Val^{sm}(V)\subset \left( (Val_i^{sm}(V)\oplus Val_{i+1}^{sm}(V)\oplus \dots
\oplus Val_n^{sm}(V)\right) \oplus Val_{i-1}^{1}(V).$$ It is
sufficient to show that $W_i(PVal_d^{sm})\cap Val_{i-1}^{1} =0$.
Assume that $\phi \in W_i(PVal_d^{sm})\cap
(Val_{i-1}^{1}(V))^{sm},\, \phi \ne 0$. By the Poincar\'e duality
 there exists $\psi \in (Val_{n-i+1}(V))^{sm}$ such that $\phi
\cdot \psi \ne 0$. We may also assume that $\psi
(K)=vol(K[n-i+1],A[i-1])$ where $A\in \ck(V)$ is fixed. By
Proposition 1.3(iii) for fixed $\psi$ the map $CVal(V)\str
CVal(V)$ given by $\xi \mapsto \xi \cdot \psi$ is a continuous map
(with the topology of uniform convergence on compact subsets of
$\ck (V)$). By using the construction of multiplication one easily
checks the following claim.
\begin{claim}
$$\psi \cdot W_i(PVal_d^{sm}) =0.$$
\end{claim}
Hence by continuity it follows that $\psi \cdot W_i(PVal_d^{sm})
=0$. We get a contradiction. \qed

\begin{proposition}
 $W_i(PVal_d^{sm})$ is  a closed subspace of $PVal_d^{sm}$.

\end{proposition}
{\bf Proof.} The statement follows from the Casselman-Wallach
theorem (Theorem \ref{casselman-wallach}) and the fact that the
target of the map $\Theta_{k,d}$ is an admissible $GL(V)$-module
of finite length. \qed


\begin{proposition}
Let $\phi \in PVal_d^{sm}$. If $\phi \in W_i( PVal_{d'}^{sm})$ for
some $d'\geq d$ then
$$\lim_{r\str +0} \frac{1}{r^{i-1}}\phi(rK+x)=0 \, \forall x\in
X,\, \forall K\in \ck (V).\eqno{(\star)}$$ Conversely if the
condition $(\star)$ holds then $\phi \in W_i(PVal_d^{sm})$.

\end{proposition}
{\bf Proof.} Let us assume that $\phi \in PVal_{d'}^{sm}$. Then
$$\int_{rK+x_0+\sum_{j=1}^{n-i}\lambda_j A_j} F(x)dx=F(x_0)
vol(rK+\sum_{j=1}^{n-i}\lambda_j A_j)+\mbox{ higher order terms
}.$$ The condition $(\star)$ follows.

Let us prove the other direction. Let $\phi \in PVal_d^{sm}$
satisfy $(\star)$. We will prove the result by induction in $d$.
For $d=0$ the statement is clear. Let us assume now that the
statement is true for valuations of degree of polynomiality less
than $d$. We have
$$\phi (K+x)=P_d(K)(x)+ \mbox{ lower order terms }.$$
Thus $P_d \in Val^{sm}\otimes Pol_d$. From the condition $(\star)$
we get:
$$0=\lim _{r\str +0} \frac{1}{r^{i-1}}\phi(rK+x)=\lim _{r\str +0}\frac{1}{r^{i-1}}
P_d(rK)(x)+\cdots .$$ Hence $P_d\in (Val_i^{sm}\oplus \dots \oplus
Val_n^{sm})\otimes Pol_d$. It follows from the Casselman-Wallach
theorem that there exists $\psi \in W_i(PVal_d^{sm})$ such that
$\psi(K+x)=P_d(K)(x)+\cdots $. Applying the assumption of
induction to the valuation $\phi -\psi$ we prove the statement.
\qed

We get immediately the following corollary.
\begin{corollary}
(i) For $d'>d$ $W_i(PVal_{d'}^{sm})\cap
PVal_d^{sm}=W_i(PVal_d^{sm}).$

(ii) Let $f:U\str V$ be an injective imbedding of linear spaces.
Then $f^*(W_i(PVal_d^{sm}(V)))\subset W_i(PVal_d^{sm}(U))$.
\end{corollary}
Using part (i) of this corollary let us define the filtration
$W_i$ on all smooth polynomial valuations $PVal^{sm}$ so that
$W_i\cap PVal_d^{sm}= W_i(PVal_d^{sm}).$

\begin{theorem}
$$W_{i_1}(PVal_{d_1}^{sm})\otimes
W_{i_2}(PVal_{d_2}^{sm})\subset W_{i_1+i_2}(PVal_{d_1
+d_2}^{sm}).$$
\end{theorem}
{\bf Proof.} It is clear from the definitions that
$$W_{i_1}(PVal_{d_1}^{sm}(V))\boxtimes
W_{i_2}(PVal_{d_2}^{sm}(W))\subset
W_{i_1+i_2}(PVal_{d_1+d_2}(V\times W)).$$ Now the statement
follows from Corollary 3.5(ii). \qed
\begin{proposition}
(i) For any $i$
$$\gamma_{i+1}\subset W_i\subset \gamma_i.$$

(ii) $W_1=\gamma_1.$
\end{proposition}
{\bf Proof.} Let us prove part (i). Let us prove first the
inclusion $\gamma_{i+1}\subset W_i$. Let $\phi\in
\gamma_{i+1}(PVal_d^{sm})$. We will prove the statement by
induction in $d$. Assume first that $d=0$. Then
$\gamma_{i+1}(Val^{sm})=(Val_i^1)^{sm} \oplus Val_{i+1}^{sm}\oplus
\dots \subset W_i(Val^{sm})=Val_i^{sm} \oplus Val_{i+1}^{sm}
\oplus \dots .$

Now let us assume that the statement holds for valuations
polynomial of degree  less than $d$. Let us prove it for $d$.
 We have
$$\phi (K+x)=P_d(K)(x)+ \mbox{ lower order terms }.$$
Thus $P_d \in \gamma_{i+1}(Val^{sm})\otimes Pol_d\subset
W_{i}(Val^{sm})\otimes Pol_d$. It follows from the case $d=0$ and
the Casselman-Wallach theorem that there exists $\psi \in
W_i(PVal_d^{sm})\cap \gamma_{i+1}(PVal_d^{sm})$ such that
$$\psi (K+x)=P_d(K)(x)+ \mbox{ lower order terms }.$$
Applying the induction assumption to the valuation $\phi -\psi$ we
obtain the result.

The second inclusion $W_i \subset \gamma_i$ follows from the fact
that if $\dim K<i$ then $\int_{K+\sum_{j=1}^{n-i}\lambda_j A_j}
\mu= O\left(\lambda^{n-i+1}\right)$ where
$\lambda:=\sqrt{\sum_{j=1}^{n-i}\lambda_j^2} \str 0$.

Let us prove part (ii). It remains to prove the inclusion $\gamma_
1\subset W_1$. Assume that $\phi \in \gamma_1\cap PVal_d^{sm}$. We
will prove by induction in $d$ that $\phi \in W_1$. For $d=0$ the
statement is clear. As previously we can write
$$\phi (K+x)=P_d(K)(x)+ \mbox{ lower order terms }$$
with $P_d \in \gamma_{1}(Val^{sm})\otimes Pol_d=
W_{1}(Val^{sm})\otimes Pol_d$. Again the Casselman-Wallach theorem
implies that there exists $\psi \in W_1(PVal_d^{sm})$ such that
$$\psi (K+x)=P_d(K)(x)+ \mbox{ lower order terms }.$$
Applying the assumption of induction to the valuation $\phi -\psi$
we prove the result. \qed

The next theorem gives an axiomatic characterization of the
filtration $W_i$.
\def\tw{\tilde W}
\begin{theorem}
Assume that we have another filtration $\{\tw_i\}$ on $PVal^{sm}$
such that

(1) $\{\tw_i\}$ is compatible with the multiplicative structure,
i.e. $\tw_i \cdot \tw_j \subset \tw_{i+j}$;

(2) $\gamma_{i+1}\subset \tw_i\subset \gamma_i \mbox{ for all }
i.$

(3) $\tw_0=\gamma_0,\,\tw_1 =\gamma_1$;

(4) $\tw_i\cap PVal_d^{sm}$ is a closed subspace of $PVal_d^{sm}$
for all $i,\, d$.

(5) $\tw_i$ is $Aff (V)$-invariant.

Then $\tw_i=W_i$.
\end{theorem}
{\bf Proof.} One has to prove that $\tw_i\cap
PVal_d^{sm}=W_i(PVal_d^{sm})$ for all $i,\, d$. As previously the
general case reduces easily to the case $d=0$. Let us prove the
statement in this case. We may assume that $i>1$. Remind that
$$\gamma_{i+1}\cap Val^{sm}=(Val_i^1)^{sm}\oplus Val_{i+1}^{sm}\oplus
\dots \oplus Val_n^{sm};$$
$$\gamma_{i}\cap Val^{sm}=(Val_{i-1}^1)^{sm}\oplus Val_i^{sm}\oplus
\dots \oplus Val_n^{sm}.$$ First let us show that $\tw_i \cap
(Val_{i-1}^1)^{sm}= 0.$ Otherwise due to the $Aff(V)$-invariance
of $\tw_i$ and $GL(V)$-irreducibility of $Val_{i-1}^1$ this
intersection will be equal to $(Val_{i-1}^1)^{sm}$. Note that it
follows from Proposition 2.2 and the Irreducibility Theorem 0.7
that $(Val_{i-1}^1)^{sm}\cdot (Val_1^1)^{sm} $ is dense in
$(Val_i^0)^{sm}$. But $(Val_1^1)^{sm}\subset \gamma_1= \tw_1$.
Also $\tw_{i}\cdot \tw_1\subset \tw_{i+1}$. Hence
$\tw_{i+1}\supset(Val_{i}^0)^{sm}$. Hence $\gamma_{i+1}\supset
(Val_{i}^0)^{sm}$. This is a contradiction.

Thus it remains to prove that $\tw_i\supset (Val_i^0)^{sm}$.
Remind that since $\tw_1=\gamma_1$ we have $(Val_1^0)^{sm}\subset
\tw_1$. Hence $((Val_1^0)^{sm})^{i}\subset \tw_i$. Again
Proposition 2.2 and the Irreducibility Theorem 0.7 imply that
$((Val_1^0)^{sm})^{i}$ is dense in $(Val_i^0)^{sm}$. So the result
follows. \qed

Let us consider the associated graded algebra
$gr_W(PVal^{sm}):=\bigoplus_{i=0}^{n}W_i/W_{i+1}.$ The next result
gives a description of it.
\begin{theorem}
There exists a canonical isomorphism of graded algebras
$$gr_W (PVal^{sm}(V))\simeq Val^{sm}(V) \otimes \CC[V]$$
where the $i$th graded term in the right side is equal to
$Val_i^{sm}(V)\otimes \CC[V]$ where $\CC[V]$ denote the algebra of
polynomial functions on $V$.
\end{theorem}
{\bf Proof.} Let  us construct the isomorphism explicitly. Let us
define a map
$$Q:W_i\str Val_i^{sm}\otimes \CC[V]$$ by $Q(\phi)(K)(x)=\lim_{r\str 0}
r^{-i} \phi(rK+x).$ First notice that $Q(\phi)$ is indeed a
translation invariant valuation. Let us check it say for $x=0$. We
have $$\phi(r(K+a))=r^iQ(\phi)(K+a)(0)+o(r^i).$$ On the other hand
$$\phi(r(K+a))=
\phi(rK+ra)=r^iQ(\phi)(K)(ra)+o(r^i)=r^iQ(\phi)(K)(0)+o(r^i).$$
Hence $Q(\phi)(K+a)(0)=Q(\phi)(K)(0)$.

Next the kernel of $Q$ is equal to $W_{i+1}$ by Proposition 3.4.
Let us show that $Q$ is surjective. Let us check that for any $d$
the map
$$Q:W_i(PVal_d^{sm})\str Val_i^{sm}\otimes \CC[V]_{\leq d}$$ is surjective
(where $\CC[V]_{\leq d}$ denote the space of polynomials of degree
at most $d$). Both spaces are admissible Fr\'echet representations
of $GL(V)$ of finite length. Hence by the Casselman-Wallach
theorem it is sufficient to show that $Q$ has dense image. Let
$\phi\in W_i(PVal_d^{sm})$ be a valuation of the form
\begin{equation} \phi(K)=\frac{\pr^{n-i}}{\pr\lam_1\dots \pr\lam_{n-i}}\bigg |_{\lam_j=0} \int_{K+\sum_{j=1}^{n-i}\lam_jA_j} f(x)dx
\end{equation} where $A_1,\dots A_{n-i}\in \ck(V)$ are fixed, and
$f$ is a polynomial of degree at most $d$. Then it is easy to see
that $Q(\phi)(K)(x)=f(x)\frac{\pr^{n-i}}{\pr\lam_1\dots
\pr\lam_{n-i}}\bigg |_{\lam_j=0} vol (K+\sum_j\lam_jA_j).$ Since
by the (proved) McMullen conjecture the valuations of the form
$\frac{\pr^{n-i}}{\pr\lam_1\dots \pr\lam_{n-i}}\bigg |_{\lam_j=0}
vol (K+\sum_j\lam_jA_j)$ are dense in $Val_i$, we deduce that $Q$
has dense image.

It remains to prove that $Q$ is a homomorphism of algebras. Let
$\phi\in W_i(PVal_{d}^{sm})$ be a representative of an element
from $W_i/W_{i+1}$ of the form (4), and $\psi \in
W_j(PVal_{d'}^{sm})$ be a representative of an element of
$W_j/W_{j+1}$ of the form $$\frac{\pr^{n-j}}{\pr\mu_1\cdots
\pr\mu_{n-j}}\bigg |_{\mu_l=0} \int_{K+\sum_{l=1}^{n-j}\mu_lB_l}
g(x)dx$$ where $B_1,\dots B_{n-j}\in \ck(V)$ are fixed, and $g$ is
a polynomial of degree at most $d'$. Then
$$(\phi\boxtimes\psi)(K)=$$
$$
\frac{\pr^{n-i}}{\pr\lam_1\cdots \pr\lam_{n-i}}
\frac{\pr^{n-j}}{\pr\mu_1\cdots \pr\mu_{n-j}}\bigg
|_0\int_{K+\sum_m\lam_m (A_m\times \{0\})+\sum_l\mu_l (\{0\}\times
B_l)} f(x)g(y)dx dy.$$ Then clearly
$$Q(\phi\boxtimes \psi)(K)((x,y))=$$
$$
f(x)g(y)\frac{\pr^{n-i}}{\pr\lam_1\cdots \pr\lam_{n-i}}
\frac{\pr^{n-j}}{\pr\mu_1\cdots \pr\mu_{n-j}}\bigg |_0
vol(K+\sum_m\lam_m (A_m\times \{0\})+\sum_l\mu_l (\{0\}\times
B_l))=$$ $$(Q(\phi)\boxtimes Q(\psi))(K)((x,y)).$$ The result
follows. \qed

\section{Further remarks.} In Theorem 2.6 of this
paper we have described the algebra
$Val^{O(n)}(\RR^n)=Val^{SO(n)}(\RR^n)$ of isometry invariant
continuous valuations on an $n$-dimensional Euclidean space
$\RR^n$. We would like to study the space of isometry invariant
continuous valuations in all dimensions simultaneously. More
precisely assume that $i_n:\RR^n \hookrightarrow \RR^{n+1}$ is the
standard isometric imbedding when the last coordinate vanishes. It
follows from the Hadwiger characterization theorem that for any
$k\leq n$ the restriction map $i_n^*:Val_k^{O(n+1)}(\RR^{n+1})\str
Val_k^{O(n)}(\RR^n)$ is an isomorphism. In other words for a fixed
$k$ the sequence $Val_k^{O(n)}\overset{i_n^*}{\leftarrow}
Val_k^{O(n+1)}\overset{i_{n+1}^*}{\leftarrow} \dots $ stabilizes.
Let us denote this limit vector space by $Val_k^{O(\infty)}$.
Consider the stable algebra of isometry invariant valuations
$$Val^{O(\infty)}:= \oplus_{k=0}^\infty Val_k^{O(\infty)}.$$
From Theorem 2.6 we easily deduce the following statement:
\begin{proposition}
The graded algebra $Val^{O(\infty)}= Val^{SO(\infty)}$ is
isomorphic to the graded algebra of polynomials in one variable
$\CC[x]$ with the grading given by the degree of a polynomial.
\end{proposition}

In \cite{alesker-univa} we have described in geometric terms the
vector space $Val^{U(m)}(\CC^m)$ of translation invariant
unitarily invariant continuous valuations on a Hermitian space
$\CC^m$. It would be of interest to describe the algebra structure
of this space (compare with Theorem 0.9). However it might be of
interest as well to describe the stable algebra of translation
invariant unitarily invariant continuous valuations
$Val^{U(\infty)}$ which is defined similarly to the previous case
using the following lemma.
\begin{lemma}
Let $i_n:\CC^n\hookrightarrow \CC^{n+1}$ be the standard Hermitian
imbedding such that the last coordinate vanishes. For any $k$ the
restriction map $i_n^*: Val_k^{U(n+1)}(\CC^{n+1})\str
Val_k^{U(n)}(\CC^n)$ is an isomorphism provided $n\geq k$.
\end{lemma}
{\bf Proof.} In \cite{alesker-gafa} we have shown that $\dim
Val_k^{U(n)}(\CC^n)=1+\min\{[k/2],n-[k/2]\}$. Hence if $n\geq k$
we have $\dim Val_k^{U(n)}(\CC^n)=\dim
Val_k^{U(n+1)}(\CC^{n+1})=1+[k/2]$. Hence it suffices to prove
that $i_n^*$ is injective. Let $\gr _k(\CC^n)$ denote the
Grassmannian of {\itshape real} $k$-dimensional subspaces of
$\CC^n$. In \cite{alesker-adv}, \cite{alesker-gafa} we have used
an imbedding of $Val_k^0(\CC^n)$ to the space $C(\gr_k(\CC^n))$ of
continuous functions on the Grassmannian as follows. Let $\phi\in
Val_k^0(\CC^n)$. For any $E\in \gr_k(\CC^n)$ consider the
restriction $\phi|_E\in Val_k(E)$. By a result by Hadwiger
\cite{hadwiger-book} $Val_k(E)$ coincides with the space of
Lebesgue measures on $E$. Hence $\phi |_E=f(E)vol_E$ where $vol_E$
is the Lebesgue measure induced by the Hermitian metric. Thus
$\phi \mapsto [E\mapsto f(E)]$ defines a map
$$Val_k^0(\CC^n)\str C(\gr_k(\CC^n)).$$
 The non-trivial fact due to D. Klain
\cite{klain-simple}, \cite{klain} is that this map is injective.
Let us consider the restriction of this map to $U(n)$-invariant
valuations. Then its image is contained in the space of
$U(n)$-invariant continuous functions on $\gr_k(\CC^n)$. Hence it
is enough to prove that the restriction map
$C(\gr_k(\CC^{n+1}))^{U(n+1)}\str C(\gr_k(\CC^n))^{U(n)}$ is
injective if $n\geq k$ where the restriction is considered under
the imbedding $\gr_k(\CC^n)\hookrightarrow\gr_k(\CC^{n+1})$. It is
enough to check that each $U(n+1)$-orbit in $\gr_k(\CC^{n+1})$
intersects non-trivially with $\gr_k(\CC^n)$. This fact follows
immediately from the explicit description of $U(n)$-orbits on
$\gr_k(\CC^n)$ in terms of K\"ahler angles due to H. Tasaki
\cite{tasaki}. \qed

We would like to notice that the next interesting case to classify
is the algebra of translation invariant continuous valuations on
the quaternionic space $\HH^n$ invariant under the group
$Sp(n)Sp(1)$.

\appendix
\section{Appendix.}
In this appendix we prove that all spaces we work with in Section
1 satisfy the assumptions of the Casselman-Wallach theorem. First
we remind the notion of representation of moderate growth
following \cite{wallach}, Ch. 11.

Let $G$ be a real reductive group. Assume that $G$ can be imbedded
into the group $GL(N,\RR)$ for some $N$ as a closed subgroup
invariant under the transposition. Let us fix such an imbedding
$p:G\hookrightarrow GL(N,\RR)$. (In our applications $G$ will be
either $GL(n,\RR)$ or a direct sum of several copies of
$GL(n,\RR)$.) Let us introduce a norm $|\cdot |$ on $G$ as
follows:

$$|g|:=\max\{p(g),p(g^{-1})\}$$
where $||\cdot||$ denotes the usual operator norm in $\RR^N$.
\begin{definition}
Let $(\pi,G,V)$ be a smooth representation of $G$ in a Fr\'echet
space $V$ (namely $V^{sm}=V$). One says that this representation
has {\itshape moderate growth} if for each continuous semi-norm
$\lambda$ on $V$ there exists a continuous semi-norm $\nu_\lambda$
on $V$ and $d_{\lambda}\in \RR$ such that
$$\lambda(\pi(g)v)\leq ||g||^{d_\lambda}\nu_{\lambda}(v)$$
for $g\in G,\, v\in V$.
\end{definition}

The proof of the next lemma can be found in \cite{wallach}, Lemmas
11.5.1 and 11.5.2.
\begin{lemma}\label{wall}
(i) If $(\pi,G,H)$ is a continuous representation of $G$ in a
Banach space $H$ then $(\pi,G,H^{sm})$ has moderate growth.

(ii) If $(\pi, G,V)$ is a representation of moderate growth. Let
$W$ be a closed $G$-invariant subspace of $V$. Then $W$ and $V/W$
have moderate growth.
\end{lemma}
The next lemma is obvious.
\begin{lemma}\label{obvious}
Let $G_1$ be a closed reductive subgroup of a reductive group $G$.
Assume that the image of $G_1$ in $GL(N,\RR)$ under the map
$p:G\hookrightarrow GL(N,\RR)$ is closed under the transposition.
Let $(\pi,G,H)$ has moderate growth. Then the restriction of this
representation to $G_1$ also has moderate growth.
\end{lemma}

In Section 1 we have discussed the space $PVal_d(V)$ of continuous
valuations on $V$ polynomial of degree $d$. Equipped with the
topology of uniform convergence on compact subsets of $\ck(V)$ it
is a Fr\'echet space. We will need a more precise statement.
\begin{lemma}\label{pval-banach}
The space $PVal_d(V)$ is a Banach space.
\end{lemma}
{\bf Proof.} Let $\phi\in PVal_d(V)$. By a result due to
Khovanskii and Pukhlikov \cite{khovanskii-pukhlikov1} the function
$t\mapsto \phi(tK)$ is a polynomial in $t\geq 0$ of degree at most
$n+d$ for any $K\in \ck(V)$. Thus
$\phi(tK)=\sum_{i=0}^{n+d}t^i\phi_i(K)$ where $\phi_i\in
PVal_d(V)$ and $\phi_i$ is $i$-homogeneous, namely
$\phi_i(tK)=t^i\phi_i(K)$ for any $\lambda
>0,\, K\in \ck(V)$. Let us denote by $PVal_{d,i}(V)$ the space of
$i$-homogeneous continuous valuations polynomial of degree at most
$d$. By the previous discussion one has
$$PVal_d(V)=\bigoplus_{i=0}^{n+d}PVal_{d,i}(V).$$
But it is clear that $PVal_{d,i}(V) $ is a Banach space with the
norm given by
$$||\phi||=\sup\{|\phi(K)| \,|\, K\in \ck(V),\, K\subset D\}$$
where $D$ is the unit ball.
\begin{proposition}\label{pval-mod}
$PVal_d^{sm}(V)$ has moderate growth as $GL(V)$-module.
\end{proposition}
{\bf Proof.} This immediately follows from Lemmas \ref{wall} and
\ref{pval-banach}. \qed
\begin{proposition}\label{prop}
In the notation of Section 1, the space $\Omega_d^n\otimes
C^\infty((\PP_+(V^*))^k,L^{\boxtimes k})$ has moderate growth as
$GL(V)$-module.
\end{proposition}

This proposition is an immediate consequence of the next more
general proposition.
\begin{proposition}\label{prop-gen}
Let the group $GL(V)$ act transitively on compact smooth manifolds
$X_1,\dots,X_k$. Let $\ce_i$ be a finite dimensional
$GL(V)$-equivariant vector bundle over $X_i$, $i=1,\dots,k$. Then
the space $C^\infty(X_1\times \dots \times X_k, \ce_1\times
\dots\times \ce_k)$ has moderate growth as $GL(V)$-module.
\end{proposition}
{\bf Proof.} Let us denote by $G:=(GL(V))^k$. The group $G$ acts
transitively on the manifold $X:=X_1\times \dots \times X_k$, and
the vector bundle $\ce:=\ce_1\times \dots \times \ce_k$ is
$G$-equivariant. Consider the Banach space $H=C(X,\ce)$. Since $G$
acts transitively on $X$ the space of $G$-smooth vectors in $H$ is
equal to $C^\infty(X,\ce)$. Hence by Lemma \ref{wall} the space
$C^\infty(X,\ce)$ has moderate growth as $G$-module. Consider the
diagonal imbedding $G_1:=GL(V)\hookrightarrow G$. By Lemma
\ref{obvious}, $C^\infty(X,\ce)$ has moderate growth as
$G_1$-module. \qed


\vskip 0.7cm

\begin{thebibliography}{99}

\bibitem{alesker-adv}
 Alesker, Semyon; On P. McMullen's conjecture on translation invariant valuations.
 Adv. Math. 155 (2000), no. 2, 239--263.
\bibitem{alesker-gafa}
Alesker, Semyon; Description of translation invariant valuations
on convex sets with solution of P. McMullen's conjecture. Geom.
Funct. Anal. 11 (2001), no. 2, 244--272.

\bibitem{alesker-icm}
Alesker, Semyon; Algebraic structures on valuations, their
properties and applications. Proceedings of ICM 2002, Beijing.

\bibitem{alesker-univa}
Alesker, Semyon; Hard Lefschetz theorem for valuations, unitarily
invariant valuations, and complex integral geometry. J.
Differential Geom., 63 (2003), 63-95. also: math.MG/0209263

\bibitem{alesker-dps}
Alesker, Semyon; On the degenerate principal series
representations of the group $GL(n,\RR)$. in preparation.

\bibitem{casselman}
Casselman, William; Canonical extensions of Harish-Chandra modules
to representations of $G$. Canad. J. Math. 41 (1989), no. 3,
385--438.

\bibitem{gelfand-vilenkin}
Gelfand, I. M.; Vilenkin, N. Ya.; Generalized functions. Vol. 4.
Applications of harmonic analysis. Translated from the Russian by
Amiel Feinstein. Academic Press [Harcourt Brace Jovanovich,
Publishers], New York-London, 1964 [1977].

\bibitem{goodey-weil}
Goodey, Paul; Weil, Wolfgang; Distributions and valuations. Proc.
London Math. Soc. (3) 49, No.3 (1984), 504-516.


\bibitem{hadwiger-book}
 Hadwiger, Hugo; Vorlesungen \"uber Inhalt, Oberfl\"ache und Isoperimetrie.
 Springer-Verlag, Berlin-G\"ottingen-Heidelberg 1957.

\bibitem{khovanskii-pukhlikov1}
Khovanskii, A. G.; Pukhlikov, A. V.; Finitely additive measures of
virtual polyhedra.
 (Russian) Algebra i Analiz 4 (1992), no. 2, 161--185;
 translation in St. Petersburg Math. J. 4 (1993), no. 2, 337--356.
 \bibitem{khovanskii-pukhlikov2}
 Khovanskii, A. G.; Pukhlikov, A. V.;
  The Riemann-Roch theorem for integrals and sums of quasipolynomials on virtual polytopes.
  (Russian) Algebra i Analiz 4 (1992), no. 4, 188--216;
 translation in St. Petersburg Math. J. 4 (1993), no. 4, 789--812.

\bibitem{klain-simple}
Klain, Daniel; A short proof of Hadwiger's characterization
theorem. Mathematika 42 (1995), no. 2, 329--339.
\bibitem{klain}
Klain, Daniel; Even valuations on convex bodies. Trans. Amer.
Math. Soc. 352 (2000), no. 1, 71--93.
\bibitem{mcmullen-euler}
 McMullen, Peter; Valuations and Euler-type relations on certain
 classes of convex polytopes.
 Proc. London Math. Soc. (3) 35 (1977), no. 1, 113--135.

\bibitem{mcmullen-conj}
 McMullen, Peter; Continuous translation-invariant valuations on
the space of compact convex sets. Arch. Math. (Basel) 34 (1980),
no. 4, 377--384.


\bibitem{mcmullen-survey}
McMullen, Peter; Valuations and dissections.
 Handbook of convex geometry, Vol. A, B, 933--988, North-Holland,
 Amsterdam, 1993.
\bibitem{mcmullen-schneider}
 McMullen, Peter; Schneider, Rolf;
Valuations on convex bodies. Convexity and its applications,
170--247, Birkh\"auser, Basel, 1983.

\bibitem{schneider-book}
 Schneider, Rolf; Convex bodies: the Brunn-Minkowski theory.
 Encyclopedia of Mathematics and its Applications, 44.
 Cambridge University Press, Cambridge, 1993.

\bibitem{tasaki}
Tasaki, Hiroyuki; Generalization of K\"ahler angle and integral
 geometry in complex projective spaces.
 Steps in differential geometry (Debrecen, 2000),
 349--361, Inst. Math. Inform., Debrecen, 2001.

\bibitem{wallach}
 Wallach, Nolan R.; Real reductive groups. I, II.
  Pure and Applied Mathematics, 132. Academic Press, Inc., Boston, MA,
  1988, 1992.
\end{thebibliography}
\end{document}